\documentclass[11pt]{amsart}

\usepackage{latexsym,amssymb, mathtools}
\usepackage{amsmath}
\usepackage{color}
\usepackage[all]{xy}

\newtheorem{Theorem}{\bf Theorem}[section]
\newtheorem{Lemma}[Theorem]{\bf Lemma} 
\newtheorem{Proposition}[Theorem]{\bf Proposition} 
\newtheorem{Corollary}[Theorem]{\bf Corollary} 
\newtheorem{Remark}[Theorem]{\bf Remark} 
\newtheorem{Example}[Theorem]{\bf Example} 
\newtheorem{Definition}[Theorem]{\bf Definition} 

\newenvironment{theorem}{\begin{Theorem}$\!\!\!$}{\end{Theorem}}
\newenvironment{lemma}{\begin{Lemma}$\!\!\!$}{\end{Lemma}}

\newenvironment{corollary}{\begin{Corollary}$\!\!\!$}{\end{Corollary}}
\newenvironment{remark}{\begin{Remark}$\!\!\!$}{\end{Remark}}

\newenvironment{definition}{\begin{Definition}$\!\!\!$}{\end{Definition}}
\newenvironment{prf}[1]
    {{\noindent \bf Proof of {#1}.}}{\hfill \qed}

\renewcommand{\qed}{\qquad\kern1pt   
\vbox{\hrule height 0.6pt      
        \hbox{\vrule width 0.6pt 
            \vbox{\vskip 6pt}  
            \hskip 3pt
            \vrule width 1.3pt} 
        \hrule depth 1.3pt}     
\kern1pt}
\allowdisplaybreaks
\numberwithin{equation}{section}

\newcommand{\dsp}{\displaystyle}   
\newcommand{\eps}{\varepsilon}          
\newcommand{\N}{{\bf N}}
\newcommand{\toinfty}{\to \infty}
\renewcommand{\l}{\left}
\renewcommand{\r}{\right}
\newcommand{\<}{\langle}
\renewcommand{\>}{\rangle}
\newcommand{\q}{\quad}
\newcommand{\qq}{\qquad}
\newcommand{\s}{\,\,\,}
\newcommand{\eq}[1]{{\begin{equation}#1\end{equation}}}

\newcommand{\eqn}[1]{\begin{equation*}#1\end{equation*}}
\newcommand{\eqsp}[1]{{\begin{equation}\begin{split}#1\end{split}\end{equation}}}

\begin{document}

\title[Characterization of completely continuous operators]
{Generalization of the Ehrling inequality \\ 
and universal characterization \\ of completely continuous operators}

\author{Mizuho Okumura}
\address[Mizuho Okumura]{Graduate School of Sciences, 
Tohoku University, Aoba, Sendai 980-8578 Japan}
\email{okumura.mizuho.p3@dc.tohoku.ac.jp}

\date{\today}

%



\maketitle

\begin{abstract}
The present work is devoted to an extension of the well-known Ehrling inequalities, which quantitatively characterize compact embeddings of function spaces, to more general operators. Firstly, a modified notion of continuity for linear operators, named \emph{Ehrling continuity} and inspired by the classical Ehrling inequality, is introduced, and then, a necessary and sufficient condition for Ehrling continuity is provided via arguments based on general topology. Secondly, general completely continuous operators between normed spaces are characterized in terms of (generalized) Ehrling type inequalities. To this end, the well-known local metrization of the weak topology (so to speak, a \emph{very weak norm}) plays a crucial role. Thanks to these results, a universal relation is observed among complete continuity, the very weak norm and generalized Ehrling type inequality. 
\end{abstract}

\section{Introduction}
We begin with recalling the definitions of 
complete continuity and compactness of 
linear operators, which 
are strictly distinguished throughout this paper.

\begin{definition}
Let $(X,\|\cdot\|_X)$ and $(Y,\|\cdot\|_Y)$ 
be normed spaces and let 
$T : X \to Y$ 
be a linear operator. 
\begin{enumerate}
    \item $T$ is said to be \emph{completely continuous} if 
    $T$ maps every weakly convergent sequence $(x_n)_n$ in $X$ 
    to a strongly convergent sequence $(Tx_n)_n$ in $Y$. 
    \item $T$ is said to be \emph{compact} if 
    $T$ maps any bounded subset of $X$ to a relatively compact  
    subset of $Y$. 
\end{enumerate}
\end{definition}

Then we make the following general remarks: 

\begin{remark}
{\rm \, 
\begin{enumerate}
\item It is well known that 
every compact operator must be completely continuous 
while the converse does not hold true generally. 
Throughout the paper, we shall explicitly distinguish 
these two notions, although they are often mixed in some literature, 
since they coincide with each other in, e.g., reflexive Banach spaces.

 \item If a linear operator $T : X \to Y$ is completely continuous, 
 then $T$ is continuous (or bounded), since every strongly convergent 
 sequence in $X$ is weakly convergent, and thus, it is mapped to 
 a strongly convergent sequence in $Y$, i.e., $T$ is continuous. 
\end{enumerate}
}
\end{remark}

In view of applications of Functional Analysis to, e.g., 
variational problems and PDEs, the complete continuity of linear operators 
often plays a crucial role, e.g., for 
constructing solutions for PDEs or for finding minimizers of functionals. 
For instance, the compactness of embeddings 
is often used to derive the strong convergence from the weak convergence 
or boundedness of a minimizing sequence or  a sequence of 
approximate solutions.

Although qualitative aspects of completely continuous operators, i.e., 
weak convergence implies strong convergence, are often discussed, 
there are less studies or attempts to reveal their \emph{quantitative} 
aspects.

Precisely speaking, we shall consider a completely continuous linear operator
\[
    T:X \to Y,
\]
where $(X,\|\cdot\|_X)$ and 
$(Y,\|\cdot\|_Y)$ are normed spaces.
Then we shall seek for a quantitative estimate for the operator $T$ such as  
\eq{\notag
\|Tu-Tv\|_Y \leqq \omega(u,v) \q \mbox{for all\s} u,v \in X,
}for some function 
$\omega : X\times X \to [0,\infty[$ 
satisfying that 
$\omega (u,v)$
is infinitesimal as  
$u \to v$ 
weakly in $X$. 
The function $\omega$ may be called  a 
\emph{modulus of continuity} for the operator $T$ and 
quantitatively measure the  uniform (complete) continuity of $T$.

\vspace{3pt}
One of well-known estimates as above for \emph{compact embeddings} is the 
so-called \emph{Ehrling inequality}, which is also known as the J.-L.~Lions lemma  (\cite[p.~173]{B}~\cite{Ehr}, \cite{S}, 
\cite[Lemma~7.6.]{Tomas}, \cite[p.~269, Theorem~3.5.]{Visintin}) and   
described as follows:

\vspace{5pt}
\noindent{\bf Ehrling inequality.}
Let $(X,\|\cdot\|_X), (Y,\|\cdot\|_Y)$ and $(Z, \|\cdot\|_Z)$ be 
normed spaces with embedding 
$X \hookrightarrow Y$ being compact and $Y \hookrightarrow Z$ being 
continuous. Then for any $\eps>0$, there exists a constant 
$C_\eps>0$ such that 
\eq{\label{ineq;Ehrling}
\|u\|_Y \leqq \eps \|u\|_X + C_\eps \|u\|_Z \qq \mbox{for all\s} u \in X. 
}

\begin{remark}
{\rm \, 
The embedding $X \hookrightarrow Y$ above is 
not just completely continuous, but 
supposed to be \emph{compact}.
}
\end{remark}

Owing to the arbitrariness of $\eps>0$ and complete continuity of 
$X \hookrightarrow Z$, 
the right-hand side of \eqref{ineq;Ehrling} 
is infinitesimal as $u \to 0$ weakly in $X$
and can be regarded as a modulus of continuity 
of the compact embedding $X \hookrightarrow Y$.

Let us extend the Ehrling inequality~\eqref{ineq;Ehrling} 
to more general linear operators, so to speak, the generalized Ehrling type 
inequality, and making use of it, let us define a new notion of continuity 
for linear operators, so to speak, the \emph{Ehrling continuity}.

\begin{definition}
Let $(X,\|\cdot\|_1), (X, \|\cdot\|_2)$ and $(Y,\|\cdot\|_Y)$
be normed spaces and 
let $T : X \to Y$ be a linear operator. 
\begin{enumerate}
    \item (generalized Ehrling type inequality) 
    The linear operator $T$ is 
    said to satisfy the 
    \emph{generalized Ehrling type inequality} 
    with the pair $(\|\cdot\|_1, \|\cdot\|_2)$, if 
    for any $\eps>0,$ 
    there exists a constant $C_\eps>0$ such that 
    \eq{\label{ineq;Ehrling type}
    \|Tu\|_Y \leqq \eps \|u\|_1 + C_\eps \|u\|_2 
    \qq \mbox{for all\s} u \in X. 
    } 
    
    \item (Ehrling continuity) The linear operator $T$ is said to 
    be \emph{Ehrling continuous} with the pair $(\|\cdot\|_1, \|\cdot\|_2)$, 
    if $T$ satisfies the generalized Ehrling type inequality~\eqref{ineq;Ehrling type} 
    with the pair
    $(\|\cdot\|_1, \|\cdot\|_2)$. 
\end{enumerate}
\end{definition}

\begin{remark}
{\rm \,
The original Ehrling inequality~\eqref{ineq;Ehrling} 
means that the compact embedding $X \hookrightarrow Y$ 
is Ehrling continuous with the pair $(\|\cdot\|_X, \|\cdot\|_Z)$, 
where we assign the norm on $X$ induced by $\|\cdot\|_Z$ 
to the (weaker) norm $\|\cdot\|_2$ (see also Corollary~\ref{cor;Ehrling} below). 
}
\end{remark}

A significance of the Ehrling type continuity (of an operator $T$)
is that, 
for any sequence $(u_n)_n$ in $X$, 
$Tu_n \to 0$ in $Y$ provided that  $\sup_{n}\|u_n\|_1<\infty$ and $\|u_n\|_2 \to 0$, 
thanks to the arbitrariness of $\eps>0$. 
In other words, 
``$\|\cdot\|_1$-boundedness" and 
``$\|\cdot\|_2$-convergence" yield
``$\|\cdot\|_Y$-convergence". 
Such a continuity of $T$ with hybrid use of two 
different norms (or topologies) in $X$ is in particular useful 
when the norm $\|\cdot\|_1$ is 
stronger and $\|\cdot\|_2$ weaker.

\vspace{5pt} 
In this paper, we shall develop a general framework to establish 
generalized Ehrling type inequalities  for completely continuous operators. 
In what follows, 
we are firstly concerned with the generalized Ehrling type 
inequality~\eqref{ineq;Ehrling type} 
from topological points of view.   
We shall obtain a new topological formulation 
which completely characterizes 
the generalized Ehrling type inequality~\eqref{ineq;Ehrling type}.
More precisely, a necessary and sufficient condition for the Ehrling continuity 
of linear operators will be provided (see Theorem~\ref{Theorem;Ehrling} below).
Then as a corollary, the classical Ehrling inequality~\eqref{ineq;Ehrling} 
is also re-proved via such a topological formulation 
(see Corollary~\ref{cor;Ehrling} below).

Secondly, the complete continuity of a linear operator
$T : X \to Y$ with $X$ having a separable dual $X^*$
is characterized in general settings in terms of the Ehrling continuity. 
Owing to the separability of the dual $X^*$, 
one can construct a norm topology $|\cdot|_\Phi$  
which induces a metric topology on $\{u \in X; \ \|u\|_X \leqq 1 \}$
equivalent to the weak topology of $X$ (see Lemma~\ref{Lemma;1} below),
where $\Phi = (\phi_k)_k$ is a dense subset of $B_{X^*}(1)$.  
Then with the help of the additional norm $|\cdot|_\Phi$, 
we shall obtain a universal characterization 
of completely continuous operators, i.e., 
a linear operator 
$T: (X, \|\cdot\|_X) \to (Y, \|\cdot\|_Y)$, with  
$(X, \|\cdot\|_X)$ having a separable dual, is 
completely continuous if and only if $T$ is 
Ehrling continuous 
with the pair $(\|\cdot\|_X, |\cdot|_\Phi)$, i.e., 
for any $\eps>0$, there exists $C_\eps>0$ such that 
\eqn{
\|Tu\|_Y \leqq \eps \|u\|_X + C_\eps |u|_\Phi \mbox{\qq for all\s} u\in X
}({see} Theorem~\ref{Theorem;1} below).

\bigskip
\noindent{\bf Notaion.}
Throughout the paper, we often use the following notation.
\begin{enumerate}
    \item For a normed space $(X,\|\cdot\|_X)$ and its dual $X^*,$ 
    the duality pairing between $X$ and $X^*$ is denoted by 
    $\<\cdot, \cdot\>.$
    
    \item For $r>0$ and a normed space $(X,\|\cdot\|_X)$,
    we denote by $B_X(r)$ the (strongly) \emph{closed} ball in $X$ of 
    radius $r$ centered at the origin, i.e., 
    $B_X(r)=\{u\in X; \ \|u\|_X \leqq r \}$.
    \item Let  $(X,d_X)$ and $(Y,d_Y)$ be metric spaces.
    We may denote  
    \eq{\notag
    T:(X,d_X) \to (Y,d_Y)  \qq \mbox{continuously},
    }in order to emphasize that $T$ is continuous 
    with respect to metrics $d_X$ and $d_Y$.
    \item For a normed space $(X,\|\cdot\|_X)$ and a subset 
    $B$ in $X$, we regard $(B,\|\cdot\|_X)$ as a metric 
    space with the metric $d$ being induced by $\|\cdot\|_X$, i.e., 
    $d(u,v)\coloneqq \|u-v\|_X$ for all $u,v \in B$. 
    \item We denote by $C$ generic non-negative constants, which do not depend on 
    elements of corresponding spaces or sets and may  vary from line to line. 
    The constant $C$ will have subscripts of parameters 
    if we emphasize the dependence of 
    $C$ on those parameters. 
\end{enumerate}


\subsection{Main results~I}
The following theorem gives a topological formulation 
which completely characterizes the generalized 
Ehrling type inequality, i.e., a necessary and sufficient condition 
for the Ehrling continuity of linear operators.

\begin{theorem}\label{Theorem;Ehrling}
Let $(X,\|\cdot\|_1), (X, \|\cdot \|_2), $ and $(Y,\|\cdot\|_Y)$ 
be normed spaces, 
let $B \coloneqq \{u\in X; \ \|u\|_1 \leqq 1 \}$  
and let $T:X \to Y$ be a linear operator. 
Assume that 
\eq{\label{theorem;Ehrling cts twist}
T|_B: (B, \|\cdot\|_2) \to (Y,\|\cdot\|_Y)
}is continuous.  
Then $T$ is Ehrling continuous with the pair $(\|\cdot\|_1, \|\cdot\|_2)$, i.e., 
for any $\eps>0$, there exists a constant $C_\eps>0$
such that 
\eq{\label{ineq;theorem Ehrling}
\|Tu\|_Y \leqq \eps \|u\|_1 + C_\eps \|u \|_2 
\qq \mbox{for all}\s u \in X.
}Conversely, if a linear operator $T$ is 
Ehrling continuous with the pair $(\|\cdot\|_1, \|\cdot\|_2)$, 
then the restricted mapping~\eqref{theorem;Ehrling cts twist}
is continuous.
\end{theorem}

\vspace{5pt}
Then we further observe that  

\begin{remark}
{\rm \, 
In addition to~\eqref{ineq;theorem Ehrling}, 
if there exists a constant $C>0$ 
such that 
$\|u\|_2 \leqq C \|u\|_1$
(respectively, $\|u\|_1 \leqq C \|u\|_2$) 
for all $u \in X$, 
then there exists a constant $C'>0$ such that 
$\|Tu\|_Y \leqq C' \|u\|_1 $ 
(respectively, $\|Tu\|_Y \leqq C' \|u\|_2$ ),  
which implies that $T$ is a bounded operator from 
$(X,\|\cdot\|_1)$ to $(Y,\|\cdot\|_Y)$
(respectively, from 
$(X,\|\cdot\|_2)$ to $(Y,\|\cdot\|_Y)$). 
}
\end{remark}

Owing to Theorem~\ref{Theorem;Ehrling}, 
one can give an alternative proof of the classical  Ehrling inequality~\eqref{ineq;Ehrling}.

\begin{corollary}
\label{cor;Ehrling} {\rm (Ehrling inequality)}
Let $(X,\|\cdot\|_X), (Y,\|\cdot\|_Y)$ and $(Z, \|\cdot\|_Z)$
be normed spaces with embeddings 
$\theta : X \hookrightarrow Y$ being compact and  
$\tau : Y \hookrightarrow Z$ being continuous. 
Set $B = \{u\in X; \ \|u\|_X \leqq 1 \}$. 
Let $X$ be also equipped with the norm $|\cdot|_X$ 
induced by $\|\cdot\|_Z$, i.e., 
$|u|_X \coloneqq \|(\tau\circ\theta)(u)\|_Z$ for all $u \in X$. 
Then the canonical injection induced by $\theta$ 
\eqn{
\theta|_B : (B, |\cdot|_X) \to (Y, \|\cdot\|_Y)
}is continuous, and therefore, 
the embedding $\theta$ is Ehrling continuous with the pair  
$(\|\cdot\|_X, |\cdot|_X)$, i.e., 
for any $\eps>0$, there exists a constant $C_\eps >0$ such that 
\eq{\label{eq; cor Ehrling}
\|\theta(u)\|_Y \leqq \eps \|u\|_X + C_\eps \|(\tau\circ\theta)(u)\|_Z
\qq \mbox{for all\s} u \in X.
}
\end{corollary}

The corollary above indicates that 
Theorem~\ref{Theorem;Ehrling} is an extension of 
the classical Ehrling lemma to more general  operators.

\begin{corollary}\label{cor of cor;Ehrling}
Under the same conditions as in Corollary~\ref{cor;Ehrling}, 
the metric spaces 
${(B, \|\theta \cdot\|_Y)}$ and $(B, |\cdot|_X)$ are homeomorphic. 
\end{corollary}

One can easily show the above corollary, since the canonical mapping 
${(B, \|\theta \cdot\|_Y)} \to (B, |\cdot|_X)$ 
is continuous by  assumption, 
and the inverse mapping 
$(B, |\cdot|_X) \to {(B, \|\theta \cdot\|_Y)} $ is also (sequentially) continuous 
thanks to the estimate~\eqref{eq; cor Ehrling}
and to the arbitrariness of $\eps>0$.

Another application of Theorem~\ref{Theorem;Ehrling} is 
concerned with 
an example given by H.~Brezis in~\cite[p.\,174, Exercise~6.13]{B} 
(hereafter we call it Brezis' example for short), 
and we shall give an alternative proof for it. 

\begin{corollary}\label{cor;Brezis} 
Let $(X,\|\cdot\|_X), (Y,\|\cdot\|_Y)$ be normed spaces and 
let $(X,\|\cdot\|_X)$ be \emph{reflexive}. 
Let $T : X \to Y$ be a linear and completely continuous operator 
and set $B = \{u \in X; \ \|u\|_X \leqq 1\}$. 
Assume that 
$X$ is also equipped with another norm $|\cdot|_X$ weaker than 
$\|\cdot\|_X$, i.e., for some $C>0,$ 
$|u|_X \leqq C\|u\|_X$ for all $u \in X.$
Then
\eqn{
T|_B : (B, |\cdot|_X) \to (Y, \|\cdot\|_Y)
}is continuous, and therefore, 
$T$ is Ehrling continuous with the pair 
$(\|\cdot\|_X, |\cdot|_X)$, i.e., 
for any $\eps>0,$ there exists a constant $C_\eps >0$ such that 
\eqn{
\|Tu\|_Y \leqq \eps \|u\|_X + C_\eps |u|_X 
\qq \mbox{for all\s} u \in X.
}
\end{corollary}

\vspace{5pt}
\begin{remark}
{\rm \, 
In Brezis' example, the reflexivity of $(X, \|\cdot\|_X)$ is crucial 
and provides a weak compactness of $B$ (see also
Lemma~\ref{Lemma:reflexive-separable} and 
the poof of 
Corollary~\ref{cor;Brezis} below), and a counter example for the 
non-reflexive case is 
given in~\cite[p.\,174, Exercise~6.13]{B}.
}
\end{remark}

\subsection{Main results~II}
We now move on to providing a general theory which can guarantee
an Ehrling type inequality for any completely continuous operator  
(whose domain has a separable dual space) by performing 
a local metrization of weak topology.

It is well known that the weak topology 
on a normed space with a separable dual
is locally metrizable (see~\cite[Theorem~3.29.]{B}
and Lemma~\ref{Lemma:weak met} below).

\vspace{5pt}

We need the following lemma on the local metrization of the weak topology.
It provides another norm having useful 
properties,
which will be employed to prove the main theorem below.

\begin{lemma}\label{Lemma;1}
Let $(X,\|\cdot\|_X)$ be a normed space 
whose dual $X^*$ is \emph{separable} and 
let $\Phi = (\phi_k)_k$ be a dense subset of $B_{X^*}(1)$. Set
$B_X(r) = \{u \in X; \ \|u\|_X \leqq r \}, \ r >0$. 
Define a norm $| \cdot |_\Phi$ on $X$ by 
\eq{\label{weak norm}
|u|_\Phi \coloneqq \sum_{k=1}^\infty 2^{-k} |\< \phi_k, u \>|
\mbox{\q for\s} u\in X. 
}Then the following holds: 
\begin{enumerate}
    \item for all $u\in X$, $|u|_\Phi \leqq \|u\|_X;$
    \item the series in the right-hand side of \eqref{weak norm} 
    converges uniformly on every bounded set
    in $(X,\|\cdot\|_X)$, i.e., 
    for any bounded subset $K$ in $(X, \|\cdot\|_X)$ and any $\eps >0$, 
    there exists a number $M = M(K, \eps) \in \N$ 
    such that
    \[
    \sum_{k=M}^\infty 2^{-k}|\<\phi_k , u\>| < \eps
    \qq\mbox{for all\s} u \in K;
    \]
    \item for a sequence $(u_n)_{n \geqq1}$ and $u$ in $X$,
    $u_n \rightharpoonup u$ weakly in $(X,\|\cdot\|_X)$ 
    if and only if 
    $\sup_{n \geqq 1} \|u_n\|_X < \infty$ and 
    $|u_n - u|_\Phi \to 0$ as $n\toinfty;$
    \item in particular, the canonical bijection
    \eqn{
    (X, \|\cdot\|_X) \to (X, |\cdot|_\Phi)
    }is completely continuous; 
    \item restricted to a closed ball $B_X(r)$ for $r>0$,  
    the metric topology of $(B_X(r),|\cdot|_\Phi)$ is equivalent to 
    the (relative) weak topology of $(B_X(r), \sigma(X,X^*));$
    \item in addition, if $X$ is  \emph{reflexive}, then 
    $(B_X(r), |\cdot|_\Phi)$ is a \emph{compact} metric space
    for $r >0$.
\end{enumerate}
\end{lemma}

We remark here that there are three notions  of convergence in  
a normed space $(X, \|\cdot\|_X)$ whose dual $X^*$ is separable: 
(i) the strong convergence in $\|\cdot\|_X$; 
(ii) the weak convergence with respect to  $\sigma(X,X^*)$; 
(iii) the additional norm convergence in $|\cdot|_\Phi$ 
defined by~\eqref{weak norm}. 
Lemma~\ref{Lemma;1} shows the following implications: 
$$
\|\cdot\|_X\mbox{-convergence} \Rightarrow
\sigma(X,X^*)\mbox{-convergence} \Rightarrow 
|\cdot|_\Phi\mbox{-convergence}.
$$
It is well known that weak convergent sequences do not always  
converge strongly in general. 
A natural question is whether or not the convergence in  
$|\cdot|_\Phi$ implies the convergence in $\sigma(X,X^*)$. 
In general, the answer is negative, 
and there is a counter example, that is, 
a sequence convergent in $|\cdot |_\Phi$ but not 
convergent in $\sigma(X,X^*)$. 
For details, see Appendix. 
Hence it is reasonable to call the norm topology $(X,|\cdot|_\Phi)$ 
a \emph{very weak topology} of $X$.

\begin{definition}\label{Definition very weak top}
{\rm (Very weak topology)}
Let $(X,\|\cdot\|_X)$ be a normed space whose dual $X^*$ 
is \emph{separable}, 
let $\Phi = (\phi_k)_k$ be a dense subset of 
$B_{X^*}(1)$ 
and let $|\cdot|_\Phi$ be 
given by~\eqref{weak norm}.
The norm topology $(X, |\cdot|_\Phi)$ is called 
the \emph{very weak topology of} $X$ 
\emph{associated with} $\Phi$ 
and the norm $|\cdot|_\Phi$ 
is called the \emph{very weak norm of}  $X$
\emph{associated with} $\Phi$.  
\end{definition}

\begin{remark}
{\rm \, 
As is seen in~\eqref{weak norm}, 
the norm $|\cdot|_\Phi$ depends on the choice of
a dense subset 
$\Phi = (\phi_k)_k$ in $B_{X^*}(1)$.
However, according to Lemma~\ref{Lemma;1}, 
every very weak norm is equivalent to the weak topology  
on every closed ball centered at the origin 
regardless of the choice of a dense subset: 
i.e.,  
the very weak topology is independent of the choice 
when restricted to a closed ball centered at the origin.  
Moreover, as will be remarked 
in Remark~\ref{Remark;theorem1} below, 
the theorem below 
will show that 
the choice makes no influence on our main result, that is, 
Ehrling type inequalities for completely continuous operators
{in terms of} the very weak norm. 
Thus it suffices to fix one dense subset of $B_{X^*}(1)$
to define and apply the very weak topology of $X$.
}
\end{remark}

Now we are in a position to describe our 
main results on a 
universal characterization of the 
complete continuity of linear operators 
defined on a normed space with a separable dual.

\begin{theorem}
\label{Theorem;1}
Let $(X, \|\cdot\|_X)$ be a normed space 
whose dual $X^*$ is \emph{separable}, 
let $(Y, \|\cdot\|_Y)$ be a normed space 
and let $\Phi=(\phi_k)_{k}$ be a dense subset of $B_{X^*}(1)$. 
Set $B_X(r) = \{u \in X; \ \|u\|_X \leqq r \}, \ r >0$. 
Define by~\eqref{weak norm} the very weak norm $|\cdot|_\Phi$ 
on $X$ associated with $\Phi$.  
Let $T:(X, \|\cdot\|_X) \to Y$ be a 
\emph{completely continuous linear operator}.
Then the following holds:
\begin{enumerate}
    \item The operator $T$ is Ehrling continuous with the pair 
    $(\|\cdot\|_X, |\cdot|_\Phi)$, i.e., 
    for any $\eps >0$, there exists a constant $C_\eps>0$ such that 
    \eq{\label{Theorem;1 interpolation 1}
    \|Tu\|_Y \leqq \eps \|u\|_X + C_\eps |u|_\Phi \qq \mbox{for all} \s u \in X.
    }
    \item If, in addition, $X$ is \emph{reflexive} and $T$ is \emph{injective}
    (in this case, 
    $X$ is \emph{separable} and \emph{reflexive} and 
    $T$ is a \emph{compact embedding}), 
    then for $r>0$, $(T(B_X(r)), \|\cdot \|_Y)$ is 
    \emph{homeomorphic} to  $(B_X(r),|\cdot|_\Phi),$ which is 
    a compact metric space, and 
    for any $\eps>0$, there exists a constant $C_\eps>0$ such that 
    for all $u \in X$, 
    \eq{\label{Theorem;1 interpolation 2}
    |u|_\Phi \leqq \eps \|u\|_X + C_\eps \|Tu\|_Y.
    }
\end{enumerate}
\end{theorem}

\vspace{5pt}
Conversely, assume that a linear operator 
$T: X \to Y$ satisfies 
the relation~\eqref{Theorem;1 interpolation 1} 
under the same settings as in Theorem~\ref{Theorem;1}. 
Then thanks to the arbitrariness of $\eps>0$, 
$T$ turns out to be completely continuous 
(weak convergence is mapped to strong convergence). 
Therefore, the relation~\eqref{Theorem;1 interpolation 1} 
gives a full characterization of 
completely continuous operators 
in terms of natural two different norms 
when $(X, \|\cdot\|_X)$ has a separable dual. 
Hence we get a quantitative estimates of completely continuous 
operators which we have sought for.

\vspace{5pt}
Moreover, we give an important remark.  

\begin{remark}\label{Remark;theorem1}
{\rm\,  
We also see that the choice of a dense subset of $B_{X^*}(1)$ 
in~\eqref{weak norm} makes no influence on the 
Ehrling continuity of linear operators. 
Let $|\cdot |_\Psi$ be the very weak norm 
associated with $\Psi=(\psi_k)_k$ which is a dense subset 
of $B_{X^*}(1)$ different from $\Phi$. 
Then the canonical bijection 
\eqn{
(X, \|\cdot\|_X) \to (X, |\cdot |_\Psi)
}is completely continuous. 
So from the above theorem, for any 
$\eps>0$, there exists a constant $C>0$ such that 
\eqsp{\notag 
|u |_\Psi &\leqq \eps \|u\|_X + C|u|_\Phi \qq \mbox{for all} \s u \in X. 
}Exchanging the roles of 
$|\cdot|_\Phi$ and $|\cdot |_\Psi$, it holds that 
for any $\eps>0,$ there is a constant $C'>0$ such that 
\eqsp{\notag 
|u|_\Phi &\leqq \eps \|u\|_X + C'|u |_\Psi \qq \mbox{for all} \s u \in X. 
}Hence generalized Ehrling type inequalities for 
completely continuous linear operators
is irrelevant to the choice of dense subsets $\Phi$ and $\Psi$ of 
$B_{X^*}(1)$
(to be precise, the constant 
$C_\eps$ in~\eqref{Theorem;1 interpolation 1} 
may depend on the choice).
}
\end{remark}

\vspace{5pt}
According to Lemma~\ref{Lemma;1}, (i) of Theorem~\ref{Theorem;1} and 
Remark~\ref{Remark;theorem1}, 
the very weak topology can be regarded as a \emph{universal} norm topology on $X$  
in the sense that 
every completely continuous operator from $X$ (to some space)
can be characterized with 
the norms $\|\cdot\|_X$ and $|\cdot|_\Phi$ 
via the generalized Ehrling type inequality, and vice versa.
Hence in a separable dual case,  we have obtained  
a universal triad: completely continuous operators,
the very weak topology and 
the generalized Ehrling type inequality.

\vspace{5pt}
Finally, 
regarding compact embeddings of reflexive 
normed spaces with separable dual spaces into normed spaces,
we shall observe
that the assertion~(ii) of Theorem~\ref{Theorem;1} provides 
more detailed information on topology than the classical Ehrling lemma 
(see also Corollary~\ref{cor of cor;Ehrling}).

\begin{corollary}
Under the same assumptions as in
the assertion~(ii) of Theorem~\ref{Theorem;1}, 
the metric spaces 
${(B_X(1), \|T \cdot\|_Y)}$ and $(B_X(1), |\cdot|_\Phi)$ 
are homeomorphic to the (relative) weak topological space 
$(B_X(1), \sigma(X,X^*))$, i.e., these three topological spaces are
homeomorphic to each other. 
\end{corollary}

Indeed, from Lemma~\ref{Lemma;1}, the metric space 
$(B_X(1), |\cdot|_\Phi)$ 
is homeomorphic to the (relative) weak 
topological space $(B_X(1), \sigma(X,X^*) )$, 
and from the estimate~\eqref{Theorem;1 interpolation 1},~\eqref{Theorem;1 interpolation 2} and the arbitrariness of $\eps>0$, 
the canonical mapping 
${(B_X(1), \|T \cdot\|_Y)} \to (B_X(1), |\cdot|_\Phi)$ and the inverse 
$(B_X(1), |\cdot|_\Phi) \to {(B_X(1), \|T \cdot\|_Y)}$ are 
both (sequentially) continuous, whence follows the topological 
equivalence.

{
In contrast to Corollary~\ref{cor of cor;Ehrling}, 
if $(X,\|\cdot\|_X)$ is a reflexive normed space with a separable dual $X^*$ and 
is compactly embedded into some normed space $(Y,\|\cdot\|_Y)$, 
it is remarkable that  the metric space 
${(B_X(1), \|T \cdot\|_Y)}$ is  
homeomorphic to $(B_X(1), \sigma(X,X^*))$  and $(B_X(1), |\cdot|_\Phi)$, 
even though the latter two spaces stand irrelevantly to 
the space $Y$. 
This observation may indicate that in this situation, 
the very weak normed space $(X, |\cdot|_\Phi)$ 
is a central (universal) space 
among normed spaces $(Y,\|\cdot\|_Y)$ 
into which $(X,\|\cdot\|_X)$ is compactly embedded, 
since 
the canonical embedding $(X,\|\cdot\|_X) \to (X, |\cdot|_\Phi)$ 
is {compact}  and 
the norm $|\cdot|_\Phi$ induces ``\emph{locally}" the same 
topology as that induced by ${\|T \cdot\|_Y}$; they are even 
equivalent to the (relative) weak topology.   
}

\vspace{5pt}
The rest of the paper is organized as follows.  
In Section~2, we shall give preliminary {facts} 
in order to prove our main results. 
In Section~3, proofs for main results are 
provided. 
In Appendix, we shall compare a  
very weak topology with strong and weak topologies: 
we will give an example 
of a sequence which converges in a very weak 
topology but does not in the strong topology.

\section{Preliminaries} \label{section:2}
In this section, we recall some preliminary facts 
in general topology and 
functional analysis on normed spaces.

\subsection{A useful lemma for homeomorphisms}
The following well-known lemma is useful to prove that 
a mapping is a homeomorphism. 

\begin{lemma}\label{lemma;homeo}
Every continuous bijection from a compact space 
onto a Hausdorff space is always a homeomorphism, i.e., 
the inverse mapping is also continuous. 
\end{lemma}

The lemma above is verified since 
such a bijection must be a closed mapping 
(i.e., closed subsets are mapped to closed subsets).

\subsection{Density argument via the Hahn-Banach theorem}
The following lemma is very useful 
in finding an element orthogonal to a strict subspace. 
For more details, we refer the reader to \cite[Corollary~1.8.]{B}.

\begin{lemma}\label{Lemma;H-B}
{\rm (\cite[Corollary~1.8.]{B})}
Let $X$ be a normed space and let $M$ be a subspace 
of $X$ such that  
$\overline{M} \neq X$. 
Then there exists $\phi \in X^*, \  \phi \neq  0,$ such that 
\eq{\notag
    \<\phi, x\>=0 \mbox{\q for all\s} x \in M.
}
\end{lemma}

\subsection{Metrizability of the weak topology}
The following lemma is one of fundamental facts 
to prove the main results of the present paper. 
For more details, we refer the reader to \cite[Theorem~3.29]{B} 
and \cite[Section~V.5.]{Dunford}.

\begin{lemma}\label{Lemma:weak met}
Let $X$ be a normed space whose topological dual 
$X^*$ is \emph{separable}. Then a closed ball 
$B_X(r)$ equipped with the (relative) weak topology $\sigma(X,X^*)$ is metrizable.
\end{lemma}

\subsection{Reflexivity and separability}
We recall relations among 
reflexivity and separability on normed spaces 
(not necessarily restricted on Banach spaces).
For proofs, we refer the reader to
\cite[p.145]{Schaefer}, \cite[pp.225--273]{Narici}
(see also \cite[Sections~3.5 and~3.6]{B} in the Banach space case).

\begin{lemma}\label{Lemma:reflexive-separable}
Let $X$ be a normed space with its topological dual $X^*$.
\begin{enumerate}
\item If $X$ is reflexive, then $X^*$ is also reflexive and 
every bounded subset of $X$ is weakly compact. 

\item $X$ is separable whenever $X^*$ is separable.
\end{enumerate} 
\end{lemma}

\section{Proofs of Theorems} 
Now we are ready to prove our main results.

\subsection{Proofs of Main results~I}
We prove Theorem~\ref{Theorem;Ehrling} and 
Corollaries~\ref{cor;Ehrling} and~\ref{cor;Brezis}.

\vspace{5pt}
\begin{prf}{Theorem~\ref{Theorem;Ehrling}}
Since the mapping $T: (B,\|\cdot\|_2) \to (Y,\|\cdot\|_Y)$ is 
continuous
particularly at the origin, 
for any $\eps>0$, there exists $\delta_\eps>0$ such that 
\eq{\label{eq;lemma2 1}
T(\{ u\in B; \  \|u\|_2 <\delta_\eps \}) \subset 
\{ u\in Y; \ \|u\|_Y < \eps \}.
}Let $u\in X\setminus \{0\}$ be arbitrary. 
One observes 
\eqn{
\frac{\delta_\eps}{\delta_\eps\|u\|_1 + \|u\|_2} u 
\in \{ u\in B; \  \|u\|_2 <\delta_\eps \}.
}Indeed,  
\eqsp{\notag 
    \frac{\delta_\eps \|u\|_1}{\delta_\eps \|u\|_1 + \|u\|_2}
    < \frac{\delta_\eps \|u\|_1}{\delta_\eps \|u\|_1} =1,  
    \qq\q
    \frac{\delta_\eps \|u\|_2}{\delta_\eps \|u\|_1 + \|u\|_2}
    < \frac{\delta_\eps \|u\|_2}{ \|u\|_2} = \delta_\eps.
}Hence due to~\eqref{eq;lemma2 1}, one has 
\eqn{
    \l\| T \l( \frac{\delta_\eps}{\delta_\eps \|u\|_1 + \|u\|_2} u \r) \r\|_Y 
    <\eps,
}which leads to 
\eq{\label{202011250010}
    \|Tu\|_Y \leqq \eps \|u\|_1 + \frac{\eps}{\delta_\eps}\|u\|_2.
}The case that  $u=0$ is trivial.

Conversely, suppose that a linear operator $T$ is Ehrling continuous with the pair 
$(\|\cdot\|_1, \|\cdot\|_2)$ and consider the restricted mapping 
\eqn{
T|_B : (B, \|\cdot\|_2) \to (Y,\|\cdot\|_Y). 
}Let $(u_n)_n$ and $u$ be in $B$  such that 
$\|u_n-u\|_2 \to 0$. From the generalized Ehrling type inequality~\eqref{Theorem;Ehrling}, 
it follows that 
\eqn{
\|Tu_n-Tu\|_Y \leqq \eps \|u_n-u\|_1 + C_\eps \|u_n-u\|_2
}and hence, passing to the limit as $n \toinfty$, one obtains 
\eqn{
\varlimsup_{n\toinfty}\|Tu_n-Tu\|_Y \leqq 2 \eps.
}Thanks to the arbitrariness of $\eps>0$, $T|_B$ turns out to be 
continuous. 
\end{prf}

\vspace{5pt}
We can obtain the classical Ehrling inequality 
as a corollary of Theorem~\ref{Theorem;Ehrling}.

\vspace{5pt}
\begin{prf}{Corollary~\ref{cor;Ehrling}}
Set $B = \{u \in X; \ \|u\|_X \leqq 1 \}$. 
We equip $X$ with the norm $|\cdot|_X$ induced by $\|\cdot\|_Z$, 
that is,  
$|u|_X \coloneqq \|(\tau\circ\theta)(u)\|_Z$ 
for all $u \in X$ where $\theta : X \to Y$ is a compact embedding and
$\tau : Y \to Z$ is a continuous embedding, and then 
we consider a metric space $(B, |\cdot |_X)$, 
for which one has a trivial homeomorphism 
\eqn{
(\tau\circ \theta)|_B : (B, |\cdot|_X) \to ((\tau\circ \theta)(B), \|\cdot\|_Z); 
\s u \mapsto  (\tau\circ\theta)(u).
}In view of Theorem~\ref{Theorem;Ehrling}, it suffices to show
that the operator 
\eqn{
\theta|_B : (B, |\cdot|_X) \to (Y, \|\cdot\|_Y); \s u \mapsto\theta (u), 
}is continuous.
To this end, we shall make use of the following commutative diagram~\eqref{CD; Ehrling}.
Symbols ``$\hookrightarrow$" and 
``$\hookrightarrow \hspace{-8pt} \rightarrow$" 
denote continuous injections and continuous bijections, respectively. 
The injections $i_1, i_2,$ and those 
with no suffix in~\eqref{CD; Ehrling}
are canonical and continuous injections. 
Also, $\dsp \underset{\cong}{\hookrightarrow \hspace{-8pt} \rightarrow}$ 
denotes a homeomorphism, i.e., a continuous bijection with a continuous inverse.

\eq{\vcenter{\label{CD; Ehrling}
\def\objectstyle{\displaystyle}
\def\objectstyle{\scriptstyle}
\def\labelstyle{\scriptstyle}
\xymatrix@M=8pt{
(X, \, \|\cdot\|_X) \ar@{^{(}->}[r]^-{\theta}_{compact} & (Y,\, \|\cdot\|_Y) \ar@{^{(}->}[r]^-{\tau} & (Z,\, \|\cdot\|_Z) &  \\
  & ( \overline{\theta(B)}^{\|\cdot\|_Y},\, \|\cdot\|_Y) \ar@{^{(}->}[u]_-{i_2} \ar@{^{(}->>}[r]^-{\tau|_{\, \overline{\theta(B)}^{\|\cdot\|_Y}}}_-{\cong}& ( \tau(\overline{\theta(B)}^{\|\cdot\|_Y}),\, \|\cdot\|_Z) \ar@{^{(}->}[u] &   \\
(B, \, \|\cdot\|_X) \ar@{^{(}->>}[r]^-{\theta|_{B}} \ar@{^{(}->}[uu] & (\theta(B),\, \|\cdot\|_Y) \ar@{^{(}->>}[r]^-{\tau|_{\theta(B)}} \ar@{^{(}->}[u] & (\tau\,\circ\,\theta(B),\, \|\cdot\|_Z) \ar@{^{(}->}[u]_-{i_1}  & %
\ar@{_{(}->>}[l]^-{\cong}_-{(\tau\,\circ\,\theta)|_B}
(B,\, |\cdot|_X)  
}
}}

Since $\theta : (X, \|\cdot\|_X) \to (Y, \|\cdot\|_Y)$ 
is compact, $\theta(B)$ is relatively compact in $Y$ which means 
$(\overline{\theta(B)}^{\|\cdot\|_Y}, \|\cdot\|_Y)$ 
is a compact metric space. 
Moreover, the restriction of $\tau$ onto $\overline{\theta(B)}^{\|\cdot\|_Y}$, 
that is, 
\eqn{
\tau|_{\, \overline{\theta(B)}^{\|\cdot\|_Y}} : 
(\overline{\theta(B)}^{\|\cdot\|_Y}, \|\cdot\|_Y) 
\to (\tau(\overline{\theta(B)}^{\|\cdot\|_Y}), \|\cdot\|_Z)
}is a continuous bijection. 
From Lemma~\ref{lemma;homeo}, we observe that  
$\tau|_{\, \overline{\theta(B)}^{\|\cdot\|_Y}}$ is a homeomorphism, 
and hence, the inverse mapping is also continuous. 
Therefore, from the commutative diagram~\eqref{CD; Ehrling}, 
one can find that 
\eqn{
\theta|_B = i_2\circ (\tau|_{\, \overline{\theta(B)}^{\|\cdot\|_Y}})^{-1}\circ i_1\circ (\tau\circ \theta)|_B : (B, |\cdot|_X) \to (Y, \|\cdot\|_Y);\s  
u \mapsto \theta (u),
}is continuous. 
Now from Theorem~\ref{Theorem;Ehrling} follows 
the conclusion.
\end{prf}

\begin{remark}
{\rm \,  
An essential point of a proof of  
the classical Ehrling lemma may be 
that $(\overline{\theta(B)}^{\|\cdot\|_Y}, \|\cdot\|_Y)$  is 
a compact metric space due to the compactness of $\theta$ and that  thus 
the mapping $\tau|_{\, \overline{\theta(B)}^{\|\cdot\|_Y}}$ 
turns out to be a homeomorphism. 
}
\end{remark}

\vspace{5pt}
\begin{prf}{Corollary~\ref{cor;Brezis}}
Set $B\coloneqq \{ u\in X; \ \|u\|_X \leqq 1\}$.
Here we note that if a linear space is furnished with 
several different norms, then there also arise different weak topologies 
corresponding to the norms. So we shall pay careful attention to 
handle such different weak topologies on $X$ and $B$. 
We denote by $\sigma_1$ and $\sigma_2$ the weak topologies 
on $X$ associated with norms $\|\cdot\|_X, |\cdot|_X$, respectively, 
and we denote by $(B, \sigma_j), \ j=1,2$,  
the topological spaces furnished with the relative topologies in $(X,\sigma_j), \ j=1,2$, 
respectively. 
We set the canonical bijection 
\eqn{
\theta : (X, \|\cdot\|_X) \to (X, |\cdot|_X); 
\s u \mapsto \theta(u)\coloneqq u. 
}Then $\theta$ is a bounded linear bijection due to the assumptions, and 
moreover, the induced linear operator 
\eqn{
\hat\theta : (X, \sigma_1) \to (X, \sigma_2); \s u \mapsto
\hat\theta (u) \equiv \theta(u),  
}is also continuous (see also \cite[Theorem~3.10.]{B}). 
We make use of the following 
commutative diagram~\eqref{CD; Brezis}.

\eq{\vcenter{\label{CD; Brezis}
\def\objectstyle{\displaystyle}
\def\objectstyle{\scriptstyle}
\def\labelstyle{\scriptstyle}
\xymatrix@M=8pt{
(Y,\, \|\cdot\|_Y)  & \ar[l]_{T}^{\substack{completely \\ continuous}} (X,\, \|\cdot\|_X) \ar@{^{(}->>}[ddd]   \ar@{^{(}->>}[rrr]^{\theta}  &  &  & (X,\, |\cdot|_X) \ar@{^{(}->>}[ddd]^{i_1} \\
&   & \ar@{_{(}->}[lu] (B,\, \|\cdot\|_X) \ar@{^{(}->>}[r]^{\theta |_{B}} \ar@{^{(}->>}[d] & (B,\, |\cdot|_X) \ar@{^{(}->>}[d]^{i_1 |_B} \ar@{^{(}->}[ru] &  \\
&   & \ar@{^{(}->}[ld]_{i_2} (B,\, \sigma_1) \ar@{^{(}->>}[r]^{\hat\theta |_{B}}_{\cong} & (B,\, \sigma_2) \ar@{_{(}->}[rd] &  \\
& \ar@{.>}[luuu]^{\hat{T}} (X,\, \sigma_1) \ar@{^{(}->>}[rrr]^{\hat\theta} &  &  & (X,\, \sigma_2) \\
}
}}All the  mappings, except for
\eqsp{\notag 
&T : (X, \|\cdot\|_X) \to (Y, \|\cdot\|_Y), \\
&\hat{T} : (X, \sigma_1) \to (Y, \|\cdot\|_Y); \s u \mapsto \hat{T}u \equiv Tu, 
}are  canonical and continuous.  
Since $T$ is completely continuous, the mapping  $\hat{T}$ defined above is 
sequentially continuous. 
It suffices to show that 
\eqn{
T|_B : (B, |\cdot|_X) \to (Y, \|\cdot\|_Y); \s u \mapsto Tu,
}is continuous. 
Since $(X, \|\cdot\|_X)$ is reflexive, 
$(B, \sigma_1)$ turns out to be compact (i.e., $B$ is weakly compact).  
Moreover,  it is well known that weak topology is  Hausdorff; 
hence so is $(B,\sigma_2)$. 
Thus Lemma~\ref{lemma;homeo} implies $\hat{\theta}|_B$ is a homeomorphism. 
Therefore, it follows that 
\eqn{
T|_B = \hat{T} \circ i_2 \circ (\hat{\theta}|_B)^{-1}\circ i_1|_B : 
(B, |\cdot|_X) \to (Y,\|\cdot\|_Y)
}is sequentially continuous between metric spaces, 
and hence, it is also continuous. 
Consequently, $T : X \to Y$ is Ehrling continuous with the pair 
$(\|\cdot\|_X,|\cdot|_X)$. 
\end{prf}

\begin{remark}
{\rm \,  
In the proof mentioned above, the following facts played crucial roles: 

\begin{enumerate}
\item  $(B, \sigma_1)$  is 
a compact space due to the reflexivity of $(X, \|\cdot\|_X)$, and thus,  
the mapping $\hat{\theta}|_B$ turns out  be a homeomorphism; 
\item  $\hat{T}$ induced by $T$ is sequentially continuous, i.e., 
$T$ is weakly-strongly sequentially continuous. 
\end{enumerate}

}
\end{remark}

\subsection{Proofs of Main results~II}
We shall prove Lemma~\ref{Lemma;1} and Theorem~\ref{Theorem;1}. 
Some of the assertions of Lemma~\ref{Lemma;1} might be well known, 
and hence, proofs for them will be omitted. 

\vspace{5pt}
\begin{prf}{Lemma~\ref{Lemma;1}}
It is well known that the nonnegative function 
$|\cdot|_\Phi$ is a norm on $X$. 
Let $K$ be a bounded subset in $(X, \|\cdot\|_X)$. It follows that 
\eqn{
|u|_\Phi \leqq \sum_{k=1}^\infty 2^{-k} \|\phi_k\|_{X^*}\|u\|_X
\leqq \|u\|_X  
\leqq \sup_{u\in K}\|u\|_X < \infty.
}Hence the series in~\eqref{weak norm} is uniformly convergent 
for  $ u \in K$. Thus the assertions~(i),~(ii) are proved. 

We next prove the assertion~(iii).
Assume that  a sequence $(u_n)$ in $X$ is weakly convergent to $u$.
Then from the Uniform Boundedness Principle, 
$$
\rho \coloneqq \sup_{n \geqq 1}\|u_n\|_X + \|u\|_X 
$$
is finite, and thus for any $\eps>0$, 
there is $M\in\N$ which depends only on $\eps$ and  $\rho$  
such that 
\eqn{
\sum_{k\geqq M+1} 2^{-k}|\<\phi_k, u_n-u\>|
\leqq \rho \sum_{k\geqq M+1} 2^{-k}
\leqq \eps \qq \mbox{for all} \s n \in \N.
}It follows that 
\eqn{
|u_n-u|_\Phi \leqq \sum_{k=1}^M 2^{-k}|\<\phi_k, u_n-u\>| + \eps,
}and hence, passing to the limit as $n\toinfty$, one gets 
\eqsp{\notag 
\varlimsup_{n\toinfty}|u_n-u|_\Phi 
&\leqq \sum_{k=1}^M 2^{-k} \varlimsup_{n\toinfty}|\<\phi_k, u_n-u\>| + \eps \\
&= \eps. 
}Thanks to the arbitrariness of $\eps >0$, one obtains 
\eqn{
\lim_{n\toinfty}|u_n-u|_\Phi =0.
}Conversely, let $\eps>0$ be arbitrary and
assume that a sequence $(u_n)_n$ and $u$ in $X$ satisfy  
$\sup_{n}\|u_n\|_X <\infty$ and $|u_n-u|_\Phi \to 0$ as $n\toinfty$. 
Fix an arbitrary $\phi \in B_{X^*}(1)$ and take $l \in \N$ 
such that  $\|\phi - \phi_l \|_{X^*} \leqq \eps.$
It follows that 
\eqsp{\notag 
|\<\phi, u_n-u\>| 
&\leqq |\<\phi-\phi_l, u_n-u\>|+ |\<\phi_l, u_n-u\>| \\
&\leqq (\sup_{n}\|u_n\|_X + \|u\|_X) \eps + 2^{l} |u_n-u|_\Phi. 
}Passing to the limit as $n\toinfty$, one gets 
\eqn{
\varlimsup_{n\toinfty}|\<\phi, u_n-u\>| 
\leqq (\sup_{n}\|u_n\|_X + \|u\|_X) \eps,
}and thanks to the arbitrariness of $\eps>0$, one obtains 
\eqn{
\lim_{n\toinfty}|\<\phi, u_n-u\>| =0.
}Thus the assertion~(iii) is verified. 

The assertion~(iv) is an immediate consequence of 
the definition of 
complete continuity and~(iii). 
The assertion~(v) is proved in~\cite[Theorem~3.29]{B}. 

As for~(vi), assume that $X$ is reflexive. 
Due to Lemma~\ref{Lemma:reflexive-separable}, 
$(B_X(r), \sigma(X,X^*))$ is a compact space, 
and thus $(B_X(r), |\cdot|_\Phi)$, homeomorphic to the above, is 
also compact.  
The proof is complete. 
\end{prf}

\vspace{5pt}

Employing Theorem~\ref{Theorem;Ehrling} and Lemma~\ref{Lemma;1}, 
we shall prove Theorem~\ref{Theorem;1}.

\vspace{5pt}
\begin{prf}{Theorem~\ref{Theorem;1}}
Set $B = \{u\in X; \ \|u\|_X \leqq 1 \}$. 
According to (v) of Lemma~\ref{Lemma;1}, the metric space 
$(B, |\cdot|_\Phi)$ is homeomorphic to $(B, \sigma(X,X^*))$, 
and $T$ is (sequentially) continuous 
from $(B, |\cdot|_\Phi)$ to $(Y,\|\cdot\|_Y)$,  
since $T$ is completely continuous. 
Hence, from Theorem~\ref{Theorem;Ehrling}, 
$T$ turns out to be Ehrling continuous with the pair $(\|\cdot\|_X, |\cdot|_\Phi)$
and the assertion~(i) is verified.

Assume that $X$ is reflexive and $T$ is injective. 
Then (vi) of Lemma~\ref{Lemma;1} implies that 
$(B_X(r), |\cdot|_\Phi)$ is a compact metric space. 
Since 
\eqn{
    T: (B_X(r), |\cdot|_\Phi) \to (Y,\|\cdot\|_Y)
}is continuous 
and every continuous mapping 
preserves compactness, one can immediately observe  that 
$(T(B_X(r)), \|\cdot\|_Y)$ is a compact metric space. 
Since the mapping 
\eqn{
    T: (B_X(r), |\cdot|_\Phi) \to (T(B_X(r)), \|\cdot\|_Y)
}is a continuous bijection, 
we can employ Lemma~\ref{lemma;homeo} 
to observe that $T$ is a homeomorphism, i.e., 
$(T(B_X(r)),\|\cdot\|_Y)$ is 
homeomorphic to $(B_X(r), |\cdot|_\Phi)$, and 
the inverse mapping 
\eqn{
    T^{-1}: (T(B_X(r)), \|\cdot\|_Y) \to (B_X(r), |\cdot|_\Phi) 
}is also continuous. 
Now for the sake of simplicity we set $r=1$ and $B=B_X(1)$,  
and then, 
for any $\eps >0$, there exists $\delta_\eps >0$ such that 
\eq{\label{eq;Theorem1 1}
T^{-1} (\{ Tu \in T(B); \  \|Tu\|_Y <\delta_\eps \})
\subset 
\{ u\in B; \  |u|_\Phi < \eps \}.
}Let $u \in X\setminus\{0\}$ be arbitrary. One has 
\eqn{
    \frac{\delta_\eps}{\delta_\eps \|u\|_X + \|Tu\|_Y} u \in \{  u\in B; \  \|Tu\|_Y <\delta_\eps \}.
}Indeed, if $u\neq 0$,  it follows that 
\eqsp{\notag 
    \frac{\delta_\eps \|u\|_X }{\delta_\eps \|u\|_X + \|Tu\|_Y} < 1, 
    \qq\q
    \frac{\delta_\eps \|Tu\|_Y}{\delta_\eps \|u\|_X + \|Tu\|_Y} < \delta_\eps.
}According to~\eqref{eq;Theorem1 1}, one obtains 
\eqn{
    \l| T^{-1} \l( \frac{\delta_\eps}{\delta_\eps \|u\|_X + \|Tu\|_Y} Tu \r) \r|_\Phi <\eps,
}which yields   
\eqn{
|u|_\Phi \leqq \eps \|u\|_X + \frac{\eps}{\delta_\eps} \|Tu\|_Y.
}The case that $u =0$ is trivial, and thus, the proof is complete.
\end{prf}


\section*{Acknowledgement}
The author would like
to thank his supervisor Goro Akagi 
for his encouragement 
and advice during the preparation of the paper.

\appendix
\section{Strict weakness of the very weak topology}
We shall prove that the very weak convergence 
(see Definition~\ref{Definition very weak top}) 
is strictly weaker than the weak convergence by employing 
Lemma~\ref{Lemma;H-B}.  
According to (iii) of Lemma~\ref{Lemma;1}, 
it suffices to show the existence of a sequence $(u_n)$ in $X$ 
such that $\|u_n\|_X \toinfty$ and $|u_n|_\Phi \to 0$. 
Such a sequence exists even if $X$ is separable and reflexive. 

\vspace{5pt}
Now we assume that $(X,\|\cdot\|_X)$ is an infinite dimensional 
reflexive and separable Banach space and we identify $X$ with $X^{**}$.
Let $(\phi_k)_k$ 
be a dense subset of $B_{X^*}(1)$ such that 
for every $k \geqq 1$, the finite dimensional subspace 
$$
M_{k} \coloneqq \mbox{span}\{\phi_1, \ldots, \phi_k  \}, 
\qq k \in \N, 
$$ 
is not dense in $X^*$, i.e., $\overline{M_k}\neq X^*$
(e.g., $X=L^2(]-1,1[)$ with the Fourier basis $(\phi_n)_n$ 
satisfies these conditions).
Applying (ii) of  Lemma~\ref{Lemma;1} with $K=B_X(1)$, 
one observes that 
for any $n\in\N$, there exists $N_n \in \N$ such that 
\eq{\label{eq;app 1}
\sum_{k \geqq N_n +1 } 2^{-k} |\< \phi_k, u  \>| < \frac{1}{n^2}
\mbox{\qq for all\s} u\in K. 
}From Lemma~\ref{Lemma;H-B}, for each $n \in \N$, 
there exists $\xi_n \in X\setminus\{0\}$ 
such that 
\eq{\label{eq;app 2}
    \<\phi, \xi_n \>=0 \mbox{\q for all\s} \phi \in M_{N_n}.
}It follows from~\eqref{eq;app 1} and~\eqref{eq;app 2} that 
\eq{\label{eq;app 3}
\l| \frac{n}{\|\xi_n\|_X} \xi_n  \r|_\Phi
=
\sum_{k =1 }^\infty 2^{-k} 
\l|\l\< \phi_k, \frac{n}{\|\xi_n\|_X} \xi_n \r\> \r| < \frac{1}{n}, 
}since $\phi_1, \ldots, \phi_{N_n} \in M_{N_n}$.
Thus one can conclude that 
\eqn{
    \l| \frac{n}{\|\xi_n\|_X} \xi_n \r|_\Phi \to 0,
    }although 
\eqn{
    \l\| \frac{n}{\|\xi_n\|_X} \xi_n \r\|_X \toinfty
}as $n\toinfty$. \qed


\end{document}